\documentclass[a4paper,10pt,DIV12]{article}
\usepackage{indentfirst}
\usepackage{latexsym}
\usepackage{amssymb, amsmath,amsfonts}
\usepackage{manfnt}
\usepackage{mathrsfs}
\usepackage[ansinew]{inputenc}

\usepackage{epsfig, color}

\newcommand{\eproof}{\hspace*{ \fill } $ \Box $ \vspace{ 0.3cm }}
\newenvironment{proof}{ \emph{Proof}. }{ \eproof }

\makeatother

\newcounter{xstep}
\newcounter{ystep}
\newcounter{xstepold}
\newcounter{ystepold}
\makeatletter
\newcommand{\MotzkinPath}[1]{
 \setcounter{xstep}{0}
 \setcounter{ystep}{0}
 \setlength{\unitlength}{5mm}
 \begin{picture}(1,1)
  \put(0,0){\circle*{0.2}}
  \@for\part:=#1
  \do{
   \setcounter{xstepold}{\value{xstep}}
   \setcounter{ystepold}{\value{ystep}}
   \stepcounter{xstep}
   \addtocounter{ystep}{\part}
   \put(\value{xstep},\value{ystep}){\circle*{0.2}}
   \put(\value{xstepold},\value{ystepold}){\line(1,\part){1}}
  }
 \end{picture}
}
\renewcommand{\SS}{\mathbb{S}}

\title{The fine structure of  the sets of involutions avoiding 4321 or 3412.}
\author{Piera Manara \thanks{Università di Parma ({\tt piera.manara@fis.unipr.it})}  and Claudio Perelli Cippo \thanks{
 Politecnico di Milano
 ({\tt claudio.perelli\_cippo@polimi.it}} .}
 \date{}
\begin{document}
\maketitle

{\textbf{Abstract}.  We study the fine structure of the sets of involutions avoiding either 4312 ($I(4321)$) or 3412 ($I(3412)$), connecting the point of view of the decomposition theorems with the one of the associated labelled Motzkin paths. The algebraic generating function of the simple involutions in $I(4321)$ is given, together with other generating functions, while the set $I(3412)$ is shown containing  no simple involutions of length $n>2$.\\
The reverse-complement bijection maintains the fine structures of $I(4321)$ and trivially of $I(3412)$.

\bigskip

\bigskip


\section{Introduction.}
Pursuing the study of the fine structure of $I(321)$ given in  \cite{1}, we analyze the fine structure of the sets of involutions avoiding 4312 ($I(4321)$),  or 3412 ($I(3412)$), connecting the use of the decomposition theorems, for which we refer to \cite{1} and \cite{AA}, with the ideas attaining to associated labelled Motzkin paths, as contained in \cite{AAA}, to which we refer for the motivation and the literature on the subject.

We recall here only the following notions.\\{ $\SS$ denotes the set of all permutations, $\SS_n$ the set of  permutations of length $n$.\\
An \textit{involution} is a permutation {$\pi$ such that $\pi(\pi(i))=i$ for all $i=1,\ldots ,n$}.\\
$I$ denotes the set of all involutions, $I_n$ the set of involutions of length $n$.}

A permutation $\pi \in \SS_{n}$ \textit{avoids the pattern  $s_{k} \in \SS_{k}$} (with $k\leq n$) if $\pi$ does not contain a subsequence order-isomorphic to $s_{k}$.

{An \textit{interval } in the permutation $\pi$ is a set of contiguous indices $\cal I$ = $[a,\,b]$, such that the set of values $\pi({\cal I }) = \left\{\pi(i):\,i \in {\cal I } \right\}$ is also contiguous. A permutation $\pi \in S_{n} $ is said to be \textit{simple } if it contains only the intervals $0,\,1,\,[1,\ldots{,}n]$.}

Given a permutation $\pi \in \SS_{n}$, the set $\left\{1,2,\,\ldots,\,n \right\}$ can be partitioned into intervals $A_{1},\ldots,A_{t}$ such that $\pi(A_{i})\,=\,A_{i},\, \forall i.$ The restrictions of $\pi$ to the intervals in the finest of these decompositions are called \textit{connected components} of $\pi.$ A permutation $\pi $ with a single connected component is called \textit{connected}.\\

Let $\sigma \in \SS_{k}$, $\alpha _1,\ldots ,\alpha _k\in\SS$: \textit{inflation} of $\sigma$ by $\alpha_1,\ldots,\alpha_k$ is the permutation $\sigma \,[\alpha_1,\ldots,\alpha_k]$ obtained by replacing each element $s_i$ of $\sigma$ by a block whose pattern is $\alpha _i$.\\ We note that permutation $\pi\in\SS$ is connected (or \textit{sum undecomposable}) if and only if it is not an inflation of 12.\\

{ An involution $\pi$ can be decomposed into disjoint cycles $$\pi\,=\,(m_1,M_1)(m_2,M_2)\ldots(m_m,M_m),$$
 whith $m_i\leq M_i$, and the $m_i$ written in increasing order. We recall that a permutation $\pi$ has an \textit{excedance} at position $i$ if $\pi(i)\,>\,i$,
a \textit{deficiency} at position $i$ if $\pi(i)\,<\,i$ and a \textit{fixed} point if $\pi (i)=i$.  Thus
we say that at the position $m_i$ there is a fixed point if $m_i=M_i$, while for $m_i<M_i$ at the positions $m_i$ there
 are the excedances (or \textit{maxima}) $M_i$, and that at the positions $M_i$ there are the deficiencies (or \textit{minima}) $m_i$. \\

A \textit{Motzkin path} of length $n$ is a lattice path starting at $(0,0)$, ending at $(n,0)$, and never going below the $x$-axis, consisting of up steps $U\,=\,(1,1)$, horizontal steps $H\,=\,(1,0)$, and down steps $D\,=\,(1,-1).$

A \textit{Dyck path} is a Motzkin path containing no horizontal steps.

An \textit{irreducible Motzkin path} is a Motzkin path that does not touch the $x$-axis except for the origin and the final destination.
 A \textit{labelling }of a Motzkin path $M$ is a map associating with every down step $D$ at height $h$ an integer $\lambda(D)$, such that $1\leq \lambda(D)\leq h$. A \textit{labelled Motzkin path} is a pair $(M,\lambda)$ where $M$ is a  Motzkin path and $\lambda$ a labelling of $M$.

 The labelling assigning the label 1 to every down step $D$ of $M$ is called \textit{unitary}; the one assigning  the label $\mu$, $\mu(D) = h(D)$ is called \textit{maximal}.

A labelled Motzkin path $(M,\lambda)$ associated with $\pi\in I$ is obtained as follows. For every $i=1,\ldots,n,$ \\
- if $i$ is a fixed point for $\pi$, take a horizontal step in the  path;\\
- if $i$ is the first element of a transposition, take an up step in the  path; \\
- if $i$ is the second element af a transposition, take a down step in the
  path, labelled with $h$, if $i$ is in the $h$-th position among integers
greater than or equal to $i$ in the cycle decomposition of $\pi.$\\

\textbf{Remark 1} Looking at the graph of an involution $\pi$, one sees
that $l+1$ is the label of a deficiency $m$ in the place $i$, with
symmetric excedance $M$,  where $l$ is the number of excedances preceding
$M$  in the graph and greater than $M$ (see figures 1 and 2).\\

In \cite{AAA} it is shown how associated labelled Motzkin paths with
unitary labelling characterize the involutions of $I(4321)$, with a
consequence in particular for $I(321)$, while the maximal labelling
characterizes $I(3412)$.}\\
The graphic interpretation given in Remark 1 allows to see very nicely
this last property
(see Appendix 2). \\

\textbf{Remark 2} An involution $\sigma$ is connected if and only if its
associated  labelled Motzkin path is irreducible.

\begin{figure} [h]
\begin{center}
  \setlength{\unitlength}{5mm}
  \begin{picture}(26,9)(0,1)

\put(0,0){ \begin{picture}(0,0)

  \put(0,3){\circle*{0.2}} \put(1,5){\circle*{0.2}} \put(2,7){\circle*{0.2}} \put(3,0){\circle*{0.2}} \put(4,4){\circle*{0.2}} \put(5,1){\circle*{0.2}} \put(6,8){\circle*{0.2}} \put(7,2){\circle*{0.2}} \put(8,6){\circle*{0.2}}
   \put(0,3){\line(1,2){2}}
   \put(2,7){\line(1,-6){1}}
   \put(3,0){\line(1,4){1}}
   \put(4,4){\line(1,-3){1}}
   \put(5,1){\line(1,6){1}}
   \put(6,8){\line(1,-6){1}}
   \put(7,2){\line(1,4){1}}
   \put(0,0){\line(1,1){8}}
   \end{picture}}
   \put(14,0){\begin{picture}(0,0)
\put(0,0){\circle*{0.2}} \put(1,1){\circle*{0.2}} \put(2,2){\circle*{0.2}} \put(3,3){\circle*{0.2}} \put(4,2){\circle*{0.2}} \put(5,2){\circle*{0.2}} \put(6,1){\circle*{0.2}} \put(7,2){\circle*{0.2}} \put(8,1){\circle*{0.2}} \put(9,0){\circle*{0.2}}
   \put(0,0){\line(1,1){3}} \put(3,3){\line(1,-1){1}}  \put(4,2){\line(1,0){1}} \put(5,2){\line(1,-1){1}} \put(6,1){\line(1,1){1}} \put(7,2){\line(1,-1){2}}
   \put(3.7,2.7){\makebox(0,0){\scriptsize $1$}}    \put(5.7,1.7){\makebox(0,0){\scriptsize $1$}}    \put(7.7,1.7){\makebox(0,0){\scriptsize $1$}}    \put(8.7,.7){\makebox(0,0){\scriptsize $1$}}
   \end{picture}
   }
  \end{picture}
 \end{center}
\caption{Plot and labelled Motzkin path of $\pi=\;468152937\;$}
\end{figure}
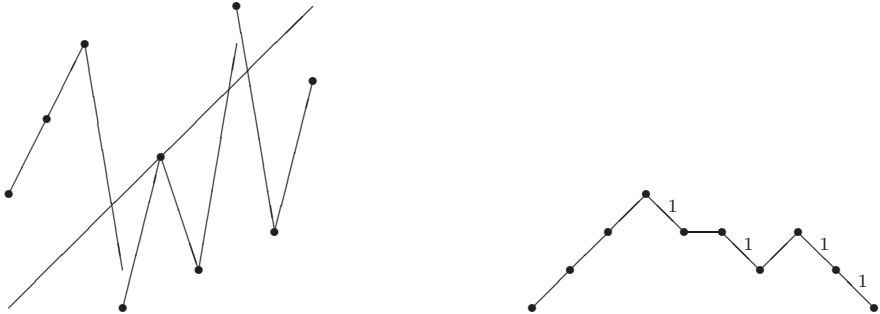

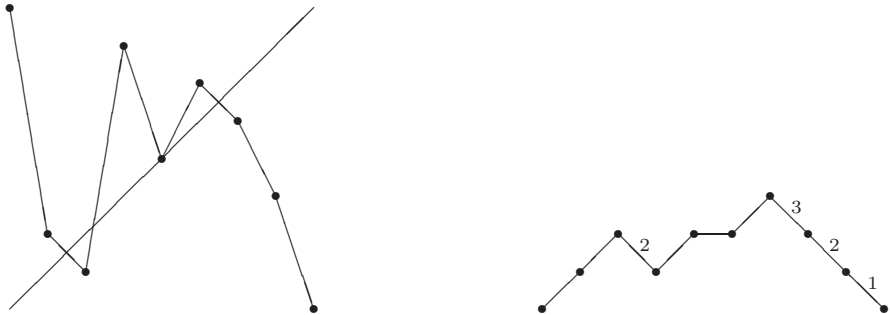
\begin{figure} [h]
\begin{center}
  \setlength{\unitlength}{5mm}
  \begin{picture}(26,9)(0,1)

  \put(0,0){\begin{picture}(0,0)
   \put(0,8){\circle*{0.2}} \put(1,2){\circle*{0.2}} \put(2,1){\circle*{0.2}} \put(3,7){\circle*{0.2}} \put(4,4){\circle*{0.2}} \put(5,6){\circle*{0.2}} \put(6,5){\circle*{0.2}} \put(7,3){\circle*{0.2}} \put(8,0){\circle*{0.2}}
   \put(0,8){\line(1,-6){1}}
   \put(1,2){\line(1,-1){1}}
   \put(2,1){\line(1,6){1}}
   \put(3,7){\line(1,-3){1}}
   \put(4,4){\line(1,2){1}}
   \put(5,6){\line(1,-1){1}}
   \put(6,5){\line(1,-2){1}}
   \put(7,3){\line(1,-3){1}}
   \put(0,0){\line(1,1){8}}
   \end{picture}}

\put(14,0){\begin{picture}(0,0)
\put(0,0){\circle*{0.2}} \put(1,1){\circle*{0.2}} \put(2,2){\circle*{0.2}} \put(3,1){\circle*{0.2}} \put(4,2){\circle*{0.2}} \put(5,2){\circle*{0.2}} \put(6,3){\circle*{0.2}} \put(7,2){\circle*{0.2}} \put(8,1){\circle*{0.2}} \put(9,0){\circle*{0.2}}
   \put(0,0){\line(1,1){2}} \put(2,2){\line(1,-1){1}} \put(3,1){\line(1,1){1}} \put(4,2){\line(1,0){1}} \put(5,2){\line(1,1){1}} \put(6,3){\line(1,-1){3}}
   \put(2.7,1.7){\makebox(0,0){\scriptsize $2$}}    \put(6.7,2.7){\makebox(0,0){\scriptsize $3$}}    \put(7.7,1.7){\makebox(0,0){\scriptsize $2$}}    \put(8.7,.7){\makebox(0,0){\scriptsize $1$}}
   \end{picture}}

  \end{picture}
 \end{center}
\caption{Plot and labelled Motzkin path of $\pi=\;932857641\;$}
\end{figure}

\section{Involutions in $I(4321)$. }

We note that  an involution in $I(321)$ is the merge of two increasing sequences of integers, the one of its  maxima (or excedances) and the one of its minima (or deficiencies), while an involution in $I(4321)$, containing $321$, is obtained by interlacing three ascending sequences of integers: the excedances, the deficiencies and the fixed points.\\

The characterization of the simple involutions in $I(321)$, exposed in \cite{1}, Theorem 8.2, can be generalized to a characterization of the  involutions in $I(4321)$.
We must briefly recall some definition we use, as in \cite{1}.\\

Consider the graph of an involution $\pi \in I(4321).$ Let's connect the points of the graph in the order their ordinates possess in the permutation $\pi$. We call \textit{plot of the involution} the drawing so obtained. \\Define two excedances (or two deficiencies) to be \textit{up-connected} (respectively \textit{down-connected}) when connected through a step of the drawing neither crossing or touching the line $y=x$ nor containing another excedance (or deficiency), therefore consecutive in the permutation. In this case, we also say that the plot has an \textit{upper connection} (or a \textit{lower connection}).\\


\noindent{\bf Proposition 2.1}{ \textit{Let $\pi \in I(4321),\,\pi \, = \, (m_1,M_1)(m_2,M_2)$ $\cdots (m_m,M_m)$.  If two  deficiencies, $m_i,\, m_{i+1}$, are down-connected, then the corresponding excedances, $M_i,\, M_{i+1}$, are consecutive integers. \\
Conversely, if two excedances, $M_i,\, M_{i+1}$, are consecutive integers, then the corresponding deficiencies, $m_i,\, m_{i+1}$, are down-connected.\\
(The same holds with excedance  replaced by deficiency, and up by down.)}

\begin{proof} Let $m_i,\, m_{i+1}$ be down-connected. If $M_i +1\neq M_{i+1}$, another integer $M$ would exist, with $M_i<M< M_{i+1}$, and the involution $\pi$ cannot have at $M$ an excedance or a fixed point, since $m_i$ and $m_{i+1}$ are down connected; $\pi$ cannot have a deficiency because $m_i$ and $m_{i+1}$ are consecutive.  \\
Conversely, if $M_i +1 = M_{i+1}$, then $m_i$ and $m_{i+1}$ are down-connected, because otherwise another excedance or fixed point $M$ should exist, with $M_i +1\,<\,M\,<\,M_{i+1}$.\end{proof}\\

{ See the involution of Figure 1, satisfying the hypothesis and the claim. On the contrary the involution of Figure 2 contains the pattern $\;4321\;$ and the claim of Proposition 2.1 is not true.}\\

The following properties of the simple involutions in $I(4321)$ are immediately obtained:
\\

 \noindent{\bf Proposition 2.2} \textit{A simple involution  $\pi \in I(4321)$ has no consecutive fixed points.}\\

\noindent{\bf Proposition 2.3} \textit{Let $\pi \in I(4321)$, $\pi$ simple. Then the plot of $\pi$ has no couples of upper and lower connections symmetric with respect to $y=x$, therefore if two excedances are up-connected, the corresponding deficiencies are not down-connected.} \\

For the inflations of 21 we have
\\\\
{ \noindent{\bf Proposition 2.4} \textit{Let $\pi \in I(4321)$ be an involution, inflation  of 21.
Then $\pi\,=\,21[\alpha_{1},\,\alpha_{2}]$, where  $\alpha_{1}=\alpha_{2}=1\,2 \ldots n$, or $\pi\,=\,321[\alpha_{1},\,\alpha_{2},\, \alpha_{3}]$, where  $\alpha_{1}=\alpha_{3}=1\,2 \ldots n$, and $\alpha_{2}=1\,2 \ldots m.$}\\
\begin{proof} It is immediate to see that otherwise $\pi$ would contain a descending sequence of length 4.
\end{proof}}

Looking at the associated Motzkin paths, the involutions in $I(4321)$ which are inflation of 21 can be obtained by expanding the involutions in $I(321)$ inflation of 21, through any number of fixed points at the maximum height, see Figure 4.\\
{ For sake of simplicity the unitary labelling will be omitted in the pictures.
The involution corresponding to a given Motkin path has as fixed points the horizontal steps and commutes every up step with the first successive down step.}\\\\

\begin{figure}[h]
\begin{center}
  \setlength{\unitlength}{4mm}
  \begin{picture}(35,5)(0,0)
   \put(0,0){\MotzkinPath{1,1,1,-1,-1,-1}}
 \put(10,0){\MotzkinPath{1,1,1,0,-1,-1,-1}}
 \put(20,0){\MotzkinPath{1,1,1,0,0,-1,-1,-1}}
 \end{picture}
 \end{center}
\caption{$\,456123\,\in I_6(321)$, $\,5674123\, \in I_7(4321)\setminus I_7(321)$, $\tau \in I(4321)\setminus I(321)$}

\end{figure}
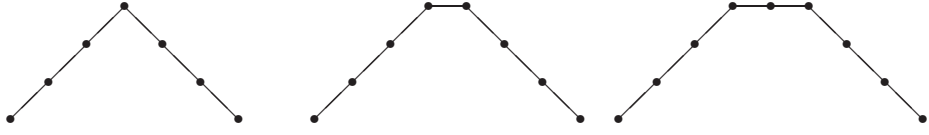

We then derive the following property:\\

\noindent{\bf Theorem 2.5} \textit{If the plot of an involution $\pi\in I(4321)$ has no couples of symmetric connections and  no consecutive fixed points, then {either} $\pi$ is  a simple involution, or $\pi$ is an involution of type 12, with $\pi\,=\,12[\alpha_{1},\,\alpha_{2}],$ where $\alpha_{1}$ is simple.}

\begin{proof} {If} $\pi$ were inflation of 21, by Proposition 2.4, its plot would have at least a couple of symmetric upper and lower connections. If $\pi$ were an {inflation } of a simple involution $\sigma \neq 12$,  its plot again would have at least a couple of symmetric upper and lower connections (deriving by the {inflation } of  a transposition, see \cite{1}, Proposition 2.6,), or at  least a couple of consecutive fixed points, or both. \\{Then either $\pi$ is simple or it is of type 12. In the last case $\pi\,=\,12[\alpha_{1},\,\alpha_{2}]$ where $\alpha_{1}$ must be simple, while $\alpha _2$, satisfying the hypothesis of the theorem, must be again either simple or of type 12}.\end{proof}

Through the use of associated labelled Motzkin paths and the characterization contained in the next Proposition, we obtain the following Theorem. \\

\noindent{\bf Proposition 2.6} \textit{(See \cite{AAA}, Theorem 3.) Let $\pi_{n}$ be an involution with $(M,\lambda)$ as the associated labelled Motzkin path of length $n$.  Then $\pi_{n}$  avoids 4321 if and only if $\lambda\,=\,\nu$ (where $\nu$ is the { unitary labelling }).}\\
As a consequence, it is shown that  \textit{$\pi_{n}$ avoids 321 if and only if $\lambda$ is the unitary labelling and all horizontal steps in $M$ are at height 0 }(see \cite{AAA}, Proposition 4).\\

\noindent{\bf Theorem 2.7}\label{T} \textit{Let $\pi_{n}\in I(4321)_{n}$ with $(M,\nu)$ as the associated labelled Motzkin path of length $n$ (with $\nu$ the unitary labelling). Then $\pi_{n}$ is simple if and only if all the following three properties hold:\\
i) $(M,\nu)$ is an irreducible  Motzkin path;\\
ii) There are no consecutive horizontal steps;\\
iii) Let $\left\{U_{1},\ldots,U_{s}\right\}$ and  $\left\{D_{1},\ldots,D_{s}\right\}$ be the sequences of the up and of the down steps in $(M,\nu)$. If two up steps $U_{i}$ and $U_{i+1}$ are consecutive up steps in $(M,\nu)$,  then the corresponding $D_{i}$ and $D_{i+1}$ are never consecutive down steps in $(M,\nu)$. }

\begin{proof} Let  $\pi_{n}\in I(4321)_{n}$ be simple, so connected: then the Motzkin path is irreducible, with no adjacent fixed points.\\
Moreover $(M,\nu)$ is such that by construction, the excedances of $\pi_n$ correspond to the up steps, the deficiencies to the down steps.  Hence, by Proposition 2.3, if $U_{i},\, U_{i+1}$ are consecutive in $(M,\nu)$, $D_{i},\, D_{i+1}$ cannot be consecutive.\\
Conversely, if $(M,\nu)$ satisfies i), ii) and iii), the involution $\pi_n\in I(4321)_{n}$ is connected, so not an inflation of 12, and the plot of $\pi_n$ has no couples of symmetric connections. It follows from Theorem 2.5 that $\pi_n$ is simple.
\end{proof}

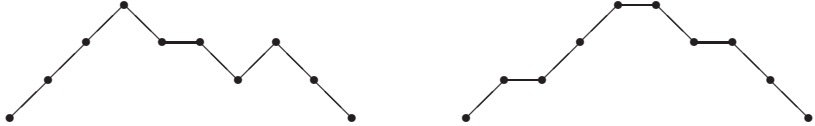
\begin{figure} [h]
\begin{center}
  \setlength{\unitlength}{4mm}
  \begin{picture}(26,5)(0,1)
   \put(0,0){\MotzkinPath{1,1,1,-1,0,-1,1,-1,-1}} \put(15,0){\MotzkinPath{1,0,1,1,0,-1,0,-1,-1}}
  \end{picture}
 \end{center}
\caption{Motzkin paths related to $\;4\,6\,8\,1\,5\,2\,9\,3\,7\;$ and to $\;6\,2\,8\,9\,5\,1\,7\,3\,4\;$}\label{fig1}

\end{figure}

The involution $\sigma =\;4\,6\,8\,1\,5\,2\,9\,3\,7\;\in I(4321)$, whose Motzkin path is illustrated in Fig.\ref{fig1}a, is simple, because the conditions i), ii), iii) are fulfilled. On the contrary $\pi =\;6\,2\,8\,9\,5\,1\,7\,3\,4\;\in I(4321)$ of Fig.\ref{fig1}b is not simple, because the pair of consecutive up steps $\;8,\,9\;$ corresponds to the pair of consecutive down steps $\;3,\,4\;$.\\

\noindent{\bf Proposition 2.8} \textit{The two following properties hold.\\
i) Starting from  a simple involution $\pi_{2n} \in I(321)_{2n}$, a simple involution $\sigma_{2n+1} \in I(4321)_{2n+1}$ can be obtained by adjoining a horizontal step to the  Dyck path $(D,\nu)$ representing $\pi_{2n}$, in any position different from the first and the last point.\\
ii) Starting from a Dyck path $(D,\nu)$ associated with an irreducible not simple involution $\pi \in I(321)$, one obtains a Motzkin path associated with a simple involution  $\sigma\in I(4321)$ by inserting in $D$ all the  horizontal steps necessary to break the consecutiveness of corresponding up and down steps.}

\begin{proof} For \textit{i)}, the thesis easily follows remembering that any simple involution, being irreducible, does  generate no horizontal steps at level zero in the associated Motzkin path, so the simple involutions of $Av(4321)$ generate horizontal steps only at levels different from zero. Adjoining a horizontal step to the  Dyck path $(D,\nu)$ (remember Proposition 2.5) associated to $\pi_{2n} \in I(321)_{2n}$ corresponds to  adjoining a fixed point, so generating a simple  involution in $Av(4321)$, as claimed.\\
While for \textit{ii)}, the thesis immediately follows from Proposition 1.6 and Theorem 1.7
\end{proof}

For instance, starting from $\pi_{6}=351624$ with associated Dyck path  that we can describe as usual in the form $UUDUDD$, if we want to insert a fixed point in the fourth position, we get the Motzkin path $UUDHUDD$, and the required involution is $3614725$.

\begin{figure} [h]
\begin{center}
  \setlength{\unitlength}{4mm}
  \begin{picture}(26,5)(0,1)
   \put(10,0){\MotzkinPath{1,1,-1,0,1,-1,-1}}
  \end{picture}
 \end{center}
\caption{Motzkin path related to $3614725$}

\end{figure}
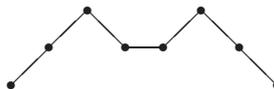

{ In \cite{AAA}, Theorem 6, (i), it is shown that the set $I_{n}(4321)$ is closed under {\it reverse-complement}; also, if the involution $\pi$ corresponds to the labelled Motzkin path $(M,\nu)$, $\nu$ the unitary labelling, the involution $\pi^{rc}$ corresponds to the path obtained by reflecting $M$ over the line $x/2$, $\nu$  unitary (see also \cite{AAA}, Proposition 2). Then, through Theorem 2.7, Proposition 2.4 and Remark 2 we obtain immediately:\\
\\
\\
{\bf Theorem 2.9} \textit{The reverse-complement bijection  preserves the fine structure of $I(4321)$.}
\\

(Let, for instance, $\sigma $ and $\pi$ be as in Fig.\ref{fig1}: then $\sigma^{rc}=371859246$, $\pi^{rc}=673951284$, as illustrated in the following picture.)
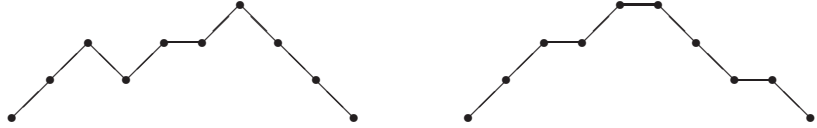
\begin{figure} [h]
\begin{center}
  \setlength{\unitlength}{4mm}
  \begin{picture}(26,5)(0,1)
   \put(0,0){\MotzkinPath{1,1,-1,1,0,1,-1,-1,-1}} \put(15,0){\MotzkinPath{1,1,0,1,0,-1,-1,0,-1}}
  \end{picture}
 \end{center}
\caption{\;$\sigma^{rc}=\;371859246\;$ and $\;\pi^{rc}=\;673951284\;$}\label{fig3}
\end{figure}
}

\section{Generating functions of subsets of I(4321).}

As we did in \cite{1}, Section 3, we consider the following generating functions of subsets of {a set of involutions:\\
$f$ denotes the generating function of the whole set;\\
$\alpha$ the generating function of the involutions which are inflation of 12;\\
$\beta$ the generating function of the involutions which are inflation of 21;\\
$\gamma$ the generating function of simple involutions different from 1, 12 e 21;\\
$\delta$ the generating function of the involutions which are {inflation } of simple involution of length  $n > 2$.} \\
It is well known that the generating function $f$ of the whole set  $I(4321)$ is $$f=-1 + \frac{1-x-\sqrt{1-2x-3x^{2}}}{2x^{2}},$$
whose expansion's coefficients are the Motzkin numbers {$1,2,4,9,21,51,127,\ldots$}\\
\\
{\bf Proposition 3.1} \textit{The generating function $\beta'$ of the involutions in $I(4321)$ inflation of 21 is $$\beta= \left(\frac{x^{2}}{1-x^{2}}\right)\left(\frac{1}{1-x}\right).$$}
\begin{proof} In fact, the thesis easily follows from the description of these involutions given in Proposition 2.4. \end{proof}
\\
\\

Always on the basis of the structure theorems  and the properties of involutions,  we then write the following relations (1) for the generating functions:\\
 $$ \begin{cases} f = x \,+\, \alpha+\beta+\gamma+\delta=-1 + \frac{1-x-\sqrt{1-2x-3x^{2}}}{2x^{2}}\\
\beta = \left(\frac{x^{2}}{1-x^{2}}\right)\left(\frac{1}{1-x}\right)\\
\alpha = (x\,+\beta+\gamma+\delta)(x\,+\,\alpha+\beta+\gamma+\delta)\end{cases}.\eqno{(1)}$$

>From (1) we obtain $$\gamma+\delta=1/4\left(2-\frac{2}{(-1+x^{2})^{2}}-\frac{3}{-1+x}-2x-\frac{1}{1+x}
-2\sqrt{1-2x-3x^{2}}\right)=$$ $$= (f+1)x^2 -\beta ,$$
whose expansion is $2x^{5}+6x^{6}+18x^{7}+47x^{8}+123x^{9}+318x^{10}+ o(x^{11}).$\\
The function $\gamma+\delta$ counts the simple involutions in $I(4321)$ and their inflations. It can obviously be obtained directly, from the following considerations.\\
The function $f=-1 + \frac{1-x-\sqrt{1-2x-3x^{2}}}{2x^{2}}= x +2x^{2}+ 4x^{3}+9x^{4}+21x^{5}+51x^{6}+127x^{7}+ \ldots$ enumerates all the involutions; by adjoining an up step and a down step respectively at the beginning and at the end of the Motzkin path one obtains a path associated to an irreducible involution; all the irreducible involutions are the ones of type $21$, the simple and their inflations. Hence the claimed property.\\

By Proposition 2.6, involutions in $I_n(4321)$ with $k$ fixed points are in bijection with Motzkin paths with $k$ horizontal steps; so the expression
 of the function $f$ in two variables can be used, enumerating the fixed points of the involutions (see \cite{AAA}, Proposition 13, and \cite{S}, sequence A097610). Precisely: $$f(x,y)=
	\frac{1-xy-\sqrt{1-2xy+x^2y^2-4x^2}}{2x^2},$$
whose expansion is	$ 1+yx+(1+y^2)x^2+(3y+y^3)x^3+(2+6y^2y^4)x^4+\ldots$,
where the coefficient $c_{k,n}$ of $c_{k,n}y^{k}x^{n}$ indicates the number of involutions of length $n$
having $k$ fixed points. (Posing $y=1$ one obtains the total number of involutions of length $n$).\\
The function $\beta$ is  immediately written in the form $$\beta(x,y)= \frac{x^{2}}{(1-x^{2})(1-xy)},$$
where the variables have the same interpretation, so giving $$(\gamma+\delta)(x,y)= \frac{1-xy-\sqrt{1-2xy+x^2y^2-4x^2}}{2}-\frac{x^{2}}{(1-x^{2})(1-xy)},$$ whose expansion is
$2yx^{5}+(1+5y^{2})x^{6}+(9y+9y^{3})x^{7}+\ldots$\\
Through the substitution $y/x$ for $y$, one obtains
$$(\gamma+\delta)(x,y)= 1/2\left(1-y-\sqrt{1-2y+y^2-4x^2}\right)-\frac{x^{2}}{(1-x^{2})(1-y)},$$ whose expansion is
$x^{6}+(2x^{4}+9x^{6})y +(5x^{4}+29x^{6})y^{2}+ (9x^{4}+69x^{6})y^{3}+\ldots,\,$
where the coefficient $c_{n,k}$ of $c_{n,k}x^{n}y^{k}$ now indicates the number of involutions of length $n+k$ having $k$ fixed points and $n/2$ transpositions.\\Using the expression of $(\gamma+\delta)(x,y)$ in the last form, one can easily obtain the function $\gamma$ enumerating the simple involutions of length greater than 2, by means of the following theorem.
\\
\\
{\bf Theorem 3.2} \label{cns} \textit{Let $\sigma \in Av(4321)$ be simple, with $k$ fixed points and $n/2$ transpositions. Then all the {inflations } of $\sigma$ belonging to $ Av(4321)$ are obtained in the two following ways:\\
\textit{i)} By inflating a fixed point with any number of new fixed points.\\
\textit{ii)} By expanding a transposition of $\sigma$ through $\alpha_{i}=\alpha^{-1}_{\sigma^{-1}(i)}=\alpha^{-1}_{\sigma(i)}\,= 123\ldots m.$}\\

\begin{proof} The assertion \textit{i)} is easily understood when thinking of the Motzkin path associated with $\sigma$.\\
As for \textit{ii)}, we observe that by inflating a transposition of $\sigma$ through a substitution presenting an inversion, we would definitely obtain an {inflation } presenting a descending sequence of length 4, so not belonging to $Av(4321)$.
\end{proof}

By means of Theorem 3.3, we can affirm that if we had the generating function  $\gamma(x,y)$ of the simple involutions in $Av(4321)$, where $x^{2}$ and $y$ count respectively the transpositions and the fixed points, we could inflate $x^{2}$ through $\frac{x^{2}}{1-x^{2}}$ and $y$ through $\frac{y}{1-y}$ in order to obtain  $(\gamma+\delta)(x,y)$, which counts  the simple involutions and their {inflations }. Because we have
$(\gamma+\delta)(x,y)$, we can go back to $\gamma(x,y)$ through the inverse functions $\frac{x^{2}}{1+x^{2}}$ and  $\frac{y}{1+y}$. We obtain in fact the following\\
\\
 {\bf Theorem 3.3}\textit{The generating function $\gamma(x,y)$ of the simple involutions in $Av(4321)$ is
$$\gamma(x,y)= 1/2\left(\frac{1}{1+y}-2x^{2}(1+y)-\sqrt{-4+\frac{4}{(1+x^{2})}+\frac{1}{(1+y)^{2}}}\right),$$ whose expansion is\\
$$(x^{6}+x^{8}+ 3x^{10}+ 6x^{12}+\ldots)+(2x^{4}+5x^{6}+ 13x^{8}+\ldots)y +(3x^{4}+14x^{6}+54x^{8}+\ldots)y^{2}+$$ $$ (x^{4}+18x^{6}+\ldots)y^{3}+ (10 x^{6}+145x^{8}+\ldots)y^{4}+\ldots$$}\\

By substituting $x$ for $y$ in $\gamma(x,y)$, we obtain the generating function in one variable:
$$\gamma(x)= 1/2\left(\frac{1}{1+x}-2x^{2}(1+x)-\sqrt{-4+\frac{4}{(1+x^{2})}+\frac{1}{(1+x)^{2}}}\right),$$
whose expansion is:\\
$2x^{5}+4x^{6}+6x^{7}+15x^{8}+31x^{9}+67x^{10}+155x^{11}+343x^{12}+787x^{13}+1829x^{14}+\ldots$\\
The sequence of the coefficients is not in Sloane \cite{S}.\\

\section{Generating functions for I(4321, 132) and for I(4321, 312). }

In \cite{AAA}, Theorem 6 and Corollary 7, the  Motzkin paths associated {with} I(4321, 132) and with I(4321, 312) are characterized and enumerative results are given. We show here how the enumerative results can be obtained through the structure theorems, and the analysis of the form of the involutions  in  $I(4321)$.\\

As usual we denote by $f$ the generating function of the whole set, by $\alpha$ the one of the involutions inflation of 12, by $\beta$ he one of the involutions inflation of 21. Moreover it is easy to see that there is no simple involution of length greater than 2 neither in $I(4321,\,132)$ nor in $I(4321,\, 312)$.
 \\Always on the basis of the structure theorems  and the { involutions' properties},  we write the following relations (2) for the generating functions of subsets of $I(4321,\,132)$:
 $$ \begin{cases} f = x \,+\, \alpha+\beta\\
\beta = \left(\frac{x^{2}}{1-x^{2}}\right)\left(\frac{1}{1-x}\right)\\
\alpha = (x\,+\beta'){x\over{1-x}}\end{cases}.\eqno{(2)}$$
In fact, any involution inflation of 21 in $I(4321)$, being of the form described in Proposition {2.4}, avoids 132, and $\beta $ is the same of Proposition 3.1.\\ The third equation says that in any involution $\pi$ of type 12, $\pi=12[\sigma_1,\sigma_2]$, the component $\sigma_2$ must be the increasing sequence $\,1\,2\,3\,\cdots\,$.  \\
>From (2) we derive $$f=-\frac{x-x^3+x^4}{(-1+x)^3(1+x)}\,=\, x+2x^2+3x^3+5x^4+\ldots ,$$
and the sequence  $1,2,3,5,7,10,13,17,...,1+\lceil{n\over 2}\rceil\,\lfloor{n\over 2}\rfloor,...$ as in the cited paper.\\\\
Now for $I(4321,\,312)$ we obtain the following relations (3):\\
$$ \begin{cases} f = x \,+\, \alpha+\beta\\
\beta = x^{2}\,+\,x^{3}\\
\alpha = (x\,+\beta)(f)\end{cases}.\eqno{(3)}$$
In fact there is no simple involution in $I(4321,\,312)$, because 1 is never a fixed point in a simple involution, so we always have a sequence of the form   $M_{i}1m_{i}$, isomorphic to 312. Moreover, by Proposition 2.4, the only involutions inflation of 21 in $I(4321)$ avoiding 312 are 21 and 321.  Finally, for any $\sigma_1,\,\sigma_2\in I(4321,321)$ the involution $\pi=12[\sigma_1,\sigma_2]$ avoids these two patterns; thus the structure theorems provide the third equation. \\
>From (3) we obtain $$f=\frac{-x-x^{2}-x^{3}}{-1+x+x^{2}+x^{3}},$$
whose expansion gives the coefficients $0,1,2,4,7,13,24,44,81...$, {(Tribonacci numbers),} as in the cited paper.\\

{\section{Generating functions for subsets of I(3412). }

For the involutions in $I(3412)$, it is well known that they are enumerated by the same function  $f=-1 + \frac{1-x-\sqrt{1-2x-3x^{2}}}{2x^{2}}$  enumerating the involutions in $I(4321)$. Through the   characterization given in \cite{AAA}, we are able to show the following properties.\\\\
{\bf Theorem 5.1} \textit{(See \cite{AAA}, Theorem 9.) Let $\pi_{n}$ be an involution with $(M,\lambda)$ as the associated labelled Motzkin path of length $n$.  Then $\pi_{n}$  avoids 3412 if and only if $\lambda$ is the maximal labelling.}\\

It is then immediate to see that the irreducible involutions in $I(3412)$ are enumerated by the same function  $(f+1)x^{2}$ enumerating the irreducible ones in $I(4321)$.\\
Also, it is easy to deduce the following useful property:\\
\\
{\bf Proposition 5.2} \textit{Let $\pi_{n}$ be an irreducible involution with $(M, \lambda)$ as the associated labelled Motzkin path of length $n$, $\lambda$ being the maximal labelling. Then one has $\pi_{n} \,= \,n \,(\pi_{n-2})\,1$.}

\begin{proof} In fact, 1 is  a deficiency of weight 1; no other deficiency may have weight 1, so the corresponding  excedance must be  $n$.
\end{proof}

It is  immediately deduced from Proposition 5.2 that every irreducible involution satisfying the hypothesis of the Proposition  is an involution inflation of  321: $\pi_{n}\,=\, 321[1,\pi_{n-2},\,1].$ So, it is easy to derive by induction the  sufficient  condition of Theorem 5.1. Moreover one can state
 the property\\\\
{\bf Theorem 5.3} \textit{ The set I(3412) has no simple involutions of length $n>2$.}
\\

{Also for Motzkin paths with maximal labelling we can omit the labelling: in this case the corresponding involution has as fixed points the horizontal steps and commutes any up step with the first successive down step {\it at the same height}.
 We can still consider the reverse-complement bijection, and in this case a claim analogous to Theorem 2.9 is trivial, because of Theorem 5.3}.\\

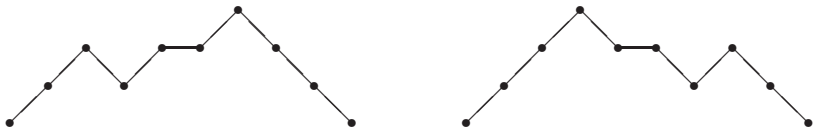
\begin{figure} [h]
\begin{center}
  \setlength{\unitlength}{4mm}
  \begin{picture}(26,5)(0,1)
   \put(0,0){\MotzkinPath{1,1,-1,1,0,1,-1,-1,-1}
   }
   \put(15,0){\MotzkinPath{1,1,1,-1,0,-1,1,-1,-1}}
  \end{picture}
 \end{center}
\caption{\;$\sigma=\;932857641\;$ and $\sigma^{rc}=\;964352871\;$ \label{fig5}}
\end{figure}

For sake of completeness, we reobtain some results given in \cite{AAA} and in Egge \cite{E}, again from the synthetic point of view of the structure theorems  and the involutions' properties. We write the following relations  for the generating functions, with the convention that $f,\alpha,\beta$ indicate the generating functions respectively of the whole considered set, of the involutions inflation of 12, of the involutions inflation of 21. Recall that $I(3412)$ has no simple involutions.\\

In $I(3412,123)$: $$ \begin{cases} f = x \,+\, \alpha+\beta\\
\beta = x^{2}\,f + x^{2}\\
\alpha = (x\,+ \frac{x^{2}}{1-x})(\frac{x}{1-x})\end{cases},\eqno{(4)}$$
obtaining the generating function $f=-\frac{x-x^3+x^4}{(-1+x)^3(1+x)}\,=\, x+x^2+3x^3+5x^4+\ldots ,$ the same of $I(4321,132)$.\\In fact, the involutions inflation of 21, avoiding the pattern 123, are precisely 21 and $321[1,\,\sigma,\,1]$, $\sigma \in Av(3412,\,123)$. While the involutions inflation of 12, avoiding the pattern 123, must be inflations of involutions showing no ascending sequences, hence the thesis.\\

In $I(3412,1234)$, we have $$ \begin{cases} f = x \,+\, \alpha+\,\beta\\
\beta = x^{2}f +x^{2}\\
\alpha = \left(-\frac{x-x^3+x^4}{(-1+x)^3(1+x)}x^{2}+x^{2} - \frac{x^{2}}{1-x}\right)\frac{x}{1-x}+\frac{x}{1-x}\frac{-x+x^3-x^4}{(-1+x)^3(1+x)}
\end{cases},\eqno{(5)}$$
obtaining the generating function $$f=-\frac{x(-1+x+x^{2}-3x^3-x^4+2x^{5}-x^{6}}{(-1+x)^5(1+x)^{2}}\,=\, x+2x^2+4x^3+8x^4+16 x^{5}+29x^{6}+51x^{7}+ \ldots $$
\\In fact, the involutions inflation of 21, avoiding the pattern 1234, are precisely 21 and $321[1,\,\sigma,\,1]$, $\sigma \in I(3412,\,1234)$. While the involutions inflation of 12, avoiding the pattern 1234, of length $n>2$, must be the following inflations:\\ a) $12[\sigma_{1},\sigma_{2}]$, where $\sigma_{1}$, involution inflation of 21, avoids 123 and contains 12, while the involution $\sigma_{2}$ avoids 12;\\
b) $12[\sigma_{1},\sigma_{2}]$, where the involution $\sigma_{1}$, inflation of 21, avoids 12, while $\sigma_{2}$ avoids 123, hence the thesis.\\

In $I(3412,132)$, we have $$ \begin{cases} f = x \,+\, \alpha+\,\beta\\
\beta = x^{2}f +x^{2}\\
\alpha= (x+\beta) \frac{x}{1-x}
\end{cases},\eqno{(6)}$$
obtaining the generating function $$f=\frac{x(1+x)}{1-x-x^2},$$ which generates the Fibonacci numbers $1,2,3,5,8, \ldots$\\
In fact, from \cite{AAA}, Theorem 10, $(i)$, one deduces that in $I(3412,132)$ the only reducible involutions of length greater than 2 are of the type $12 [\sigma_1,\,\sigma_2]$ where $\sigma_2$ is the increasing sequence $\sigma_2=\,1\,2\,\cdots\,$, hence the thesis. \\Analogously for $I(3412,213)$, the only reducible involutions of length greater than 2 have $\sigma_1=\,1\,2\,\cdots\,$, hence again one easily obtains the same equations (6). }

\section{Appendix. Patterns for the simple involutions of I(4321), for $n = 5,6,7,8,9,10.$ }

{Using $\gamma'(x,y)$, remembering Proposition 7 we give the following description for the patterns of the simple involutions of I(4321).} \\

\textit{For $n=5$}: one has 2 involutions with 1 fixed point, obtained through $\beta_{4}=3412 \in Av(321)${, expanding $\beta$ through a fixed point in  such a way as to break the consecutiveness of up steps and corresponding down  steps: }\\
4\textbf{2}513\\
351\textbf{4}2\\

\textit{For $n=6$}: there are the involution $\sigma_{6} = 351624 \in Av(321)$ and 3 involutions with 2 fixed points, obtained through $\beta_{4}\in Av(321):$\\
4\textbf{2}61\textbf{5}3\\
5\textbf{2}6\textbf{4}13\\
46\textbf{3}1\textbf{5}2\\

\textit{For $n=7$}: there are 5 involutions with 1 fixed point, obtained through 351624 {inserting a single fixed point, (remember Theorem 6)}, and 1 involution with 3 fixed points, obtained through $\beta_{4}\in Av(321):$\\
4\textbf{2}61735\\
46\textbf{3}1725\\
361\textbf{4}725\\
3617\textbf{5}24\\
35172\textbf{6}4, and\\
5\textbf{2}7\textbf{4}1\textbf{6}3\\

\textit{For $n=8$}: one has $\sigma_{8} = 35172846 \in Av(321)$ and 14 involutions with 2 fixed points, of whom 10 obtained through 351624, {inserting 2 fixed points,}  and 4 obtained through $\beta_{6}=\,456123 \in Av(321)${, always inserting the fixed points in such a way as to break the consecutiveness: for example}\\
5\textbf{2}7\textbf{4}1836\\
.........\\
6\textbf{2}7\textbf{4}8135\\

\textit{For $n=9$}: \textit{13 involutions} with 1 fixed point, of whom 7 obtained through the simple 35172846, and 6 through the 3 {inflations } of 351624 of length 8 (2 for each inflation), the insertion of a fixed point breaking the consecutiveness: for example\\
4\textbf{2}6183957\\
..........\\
while through the {inflation } 45712836 we obtain
5\textbf{2}6813947 and 4681\textbf{5}2937\\
.........\\

\textit{18 involutions} with 3 fixed points, of whom 10 obtained through 351624, and 8 through $\beta_{6}\in Av(321):$ for example\\
5 \textbf{2}7\textbf{4}193\textbf{8}6\\
..........\\
6 \textbf{2}8\textbf{4}91\textbf{7}35 and 6\textbf{2 }7\textbf{4}913 \textbf{8}5\\
..........\\

\textit{For $n = 10$}: one has the 3 simple involutions belonging to $Av(321)$;\\

\textit{54 involutions} with 2 fixed points, of whom 21 obtained through 35172846, and 33 through the 3 {inflations } of 351624 of length 8 (11 for each {inflation}): for example\\
5274193(10)68\\
..............\\
while through the {inflation} 45712836 we obtain \\
6274913(10)58, 6279513(10)48,.............................\\

\textit{10 involutions} with 4 fixed points, of whom 5 obtained through 351624, and 5 through  $\beta_{6}\in Av(321):$ for example\\
529416(10)837\\
..............\\
6284(10)17395\\
...............\\
\section{Appendix 2. An easy graphic interpretation of Theorem 5.1. }

We want to show that $\pi \in I(3412)$ has the associated Motzkin path
$(M,\lambda)$ with $\lambda$ maximal label.\\
We recall that the \textit{height} $h(D)$ of a down step $D$, associated
with a deficiency $m$, is the difference between the number of excedances
and the number of deficiencies preceding $m$ in the permutation.

If $(M,\lambda)$ had a down step $D$, corresponding to the deficiency $m$
of symmetric excedance $M$, with $\lambda(D) < h(D)$, then there would be
in $\pi$ an excedance $M_{1}$ necessarily following $M$ in the permutation
and preceding $m$, with $M \,< \,M_{1}$. The deficiency $m_{1}$ then
necessarily follows $m$, with $m_{1}\,<\, M$, otherwise $m$ would precede
$M_{1}$ in the permutation. We so deduce the sequence
$M\,M_{1}\,m\,m_{1}$, isomorphic to 3412.

Conversely, let the involution $\pi$ contain a pattern $\,abcd\,$ isomorphic to $\,3412\,$, so satisfying the inequalities $$\begin{cases} c<d<a<b\, ,\\\pi^{-1}(a)<\pi^{-1}(b)<\pi^{-1}(c)<\pi^{-1}(d),\end{cases}$$
from which we can always obtain a new isomorphic pattern of the form \\
$\,M_1M_2m_1m_2\,$ with $(m_1,M_1)\,,\,(m_2,M_2)\,$ transpositions. The deficiency $m_1$ does not have maximum label because its symmetric $M_1$ is followed by the excedance $M_2$, whose deficiency $m_2$ follows $m_1$.

\end{document}